\documentclass[12pt]{article}
\usepackage{amsmath}
\usepackage{latexsym,amssymb}
\usepackage{amsmath,amsthm}
\usepackage{pstricks}
\usepackage{delarray}
\usepackage{tabularx}
\usepackage{amsbsy}
\usepackage{graphicx}
\usepackage{amsfonts}
\usepackage{layout}
\usepackage{mathrsfs}

\def\no{\noindent}
\input cyracc.def

\begin{document}
\title{\bf N-derivations for finitely generated graded Lie algebras
}
\author{ Cui Chen\\
{\footnotesize School of Mathematics and Physics, Fujian  University
of Technology, Fuzhou 350108,
China}\\
\vspace{5mm} {\small  Email: chencui@fjut.edu.cn}\\
{ Haifeng  Lian}\\
{\footnotesize Mathematics Department, Fujian Agriculture and  Forestry University, Fuzhou 350002, China}\\
 {\small Email: lianhaif@qq.com}}
\date{}
\maketitle \no{\bf Abstract:} $N$-derivation is the natural
generalization of derivation and  triple derivation. Let ${\cal L}$
be a finitely generated Lie algebra graded by a  finite dimensional
Cartan subalgebra. In this paper, a sufficient condition for Lie
$N$-derivation algebra of ${\cal L}$ coinciding with Lie derivation
algebra of ${\cal L}$ is given. As applications, any $N$-derivation
of Schr\"{o}dinger-Virasoro algebra, generalized Witt algebras,
Kac-Moody algebras and their Borel subalgebras, is a derivation.

\vspace{3mm} \no{\bf Keywords:}\quad
 Derivation, Cartan subalgebra, Virasoro algebra, Kac-Moody algebra

\vspace{3mm} \noindent{\bf 2000 MR Subject Classification}\quad
17B40, 17B67, 17B68

 \section{Introduction}

    Let $\cal L$ be a Lie algebra over an arbitrary field $F$. Recall that an $F$-linear mapping $\psi : \mathcal{L}\rightarrow
 \mathcal{L}$ is called a Lie derivation of $\mathcal{L}$ if
$$\psi([x,y])=[\psi(x),y]+[x,\psi(y)],\hspace{1cm} \forall x,y\in
\mathcal{L}.
$$

Let $N\geq2$ be a positive integer.
 We say $\varphi$ is an
\emph{$N$-derivation} of ${\cal L}$, if $\varphi$ is a linear map
from ${\cal L}$ to itself, and satisfies
$$\varphi([x_1, \cdots, x_{N-1}, x_N])=\sum_{i=1}^{N}[x_1,\cdots,x_{i-1},\varphi(x_{i}),x_{i+1},\cdots,
x_N]\quad (\forall x_1,\cdots, x_N\in {\cal L}),$$ where $[x_1,
\cdots, x_{N-1}, x_N]=[x_1 \cdots [x_{N-2},[x_{N-1},x_N]]\cdots]$.

The set of Lie $N$-derivation is clearly a Lie algebra under the
usual bracket and will be denoted by
$\textrm{Der}^{(N)}\mathcal{L}$. By the definition, we have
$\textrm{Der}^{(2)}\mathcal{L}=\textrm{Der}\mathcal{L}$, where
$\textrm{Der}\mathcal{L}$ is the Lie derivation algebra of
$\mathcal{L}$. A Lie derivation is obviously a Lie $N$-derivation,
which implies
$\textrm{Der}\mathcal{L}\subseteq\textrm{Der}^{(N)}\mathcal{L}$. In
general, the derivation algebra $\textrm{Der}\mathcal{L}$ is a
proper subalgebra of $\textrm{Der}^{(N)}\mathcal{L}$ for $N\geq 3$.
It would be interesting to know when these two algebras coincide.

For $N=3$, Lie $N$-derivation is called Lie triple derivation which
has received a fair amount of attentions (see
\cite{JW,L,M,WL,WY,ZLF}). In \cite{WY},  Wang and Yu showed that a
linear map on Borel subalgebra of finite-dimensional simple Lie
algebra over an algebraically closed field  $F$ of characteristic
zero is a Lie triple derivation, if and only if it is an inner
derivation. Which are examples for Lie triple derivation algebra
coinciding with Lie derivation algebra, and the Lie algebras
concerned are  graded by a finite-dimensional Cartan subalgebra.

The Schr\"{o}dinger-Virasoro algebra $\mathfrak{sv}$  introduced by
Henkel in \cite{H}, during his study on the invariance of the free
Schr\"{o}dinger equation, is an infinite-dimensional Lie algebra
with $\mathbb{C}$-basis $\{L_n,M_n,Y_{n+\frac{1}{2}},C \mid n\in
\mathbb{Z}\}$ subject to the following Lie brackets:
\begin{align*}
&[L_m,L_n]=(n-m)L_{n+m}+\delta_{m+n,0}\frac{n^3-n}{12}C,\quad
[L_m,M_n]=nM_{n+m},\\
&[L_m,Y_{n+\frac{1}{2}}]=(n+\frac{1-m}{2})Y_{m+n+\frac{1}{2}},\quad [Y_{m+\frac{1}{2}},Y_{n+\frac{1}{2}}]=(n-m)M_{m+n+1},\\
&[M_m,M_n]=[M_m,Y_{n+\frac{1}{2}}]=0,\quad [\mathfrak{sv},C]=\{0\}.
\end{align*}
Due to its important applications in many areas of Mathematics and
Physics, the structure and representation theory of  $\mathfrak{sv}$
have been extensively studied (see \cite{GJP,LY1,RU,U}. For example,
in \cite{RU}, Rosen and Unterberger presented detailed cohomological
study and determined that $\mathfrak{sv}$ has three linear
independent  outer derivations. This follows that not all derivation
of $\mathfrak{sv}$ is an inner derivation.

In this paper, we consider Lie $N$-derivations for finitely
generated graded Lie algebra  ${\cal L}$. In section 2, we prove
that the Lie $N$-derivation algebra of ${\cal L}$ is graded, which
generalizes the result obtained by Farnsteiner (see \cite{Fa}). In
section 3, we assume  ${\cal L}$ is graded by a finite-dimensional
Cartan subalgebra and is over an algebraically closed field $F$ of
characteristic zero. We show that
$\textrm{Der}^{(N)}\mathcal{L}=\textrm{Der}\mathcal{L}$ for $N\geq
3$, if ${\cal L}$ satisfies condition (P) (see theorem 3.2). As
applications,  in section 4, we show that the Lie $N$-derivation
algebras of Schr\"{o}dinger-Virasoro algebra,  generalized Witt
algebras and Kac-Moody algebras are coinciding with the derivation
algebras. In particularly, any $N$-derivation of finite-dimensional
semisimple complex Lie algebra is an inner derivation.

\section{$N$-derivation of finitely
generated graded Lie algebra}

 Let $G$ be an
abelian group, and let ${\cal L}=\oplus _{{\bf \alpha}\in G}{\cal
L}_{{\bf \alpha}}$ be a $G$-graded Lie
 algebra. For $x_1,\cdots,x_n\in {\cal L}$, set $$[x_1,
\cdots, x_{n-1}, x_n]=[x_1 \cdots [x_{n-2},[x_{n-1},x_n]]\cdots].$$

Let $N\geq2$ be a positive integer. As above, denote
$\textrm{Der}^{(N)}({\cal L})$, $\textrm{Der}({\cal L})$ the set of
all $N$-derivations of ${\cal L}$  and  the set of all derivations
of ${\cal L}$, respectivly. We have $\textrm{Der}({\cal
L})\subseteq\textrm{Der}^{(N)}({\cal L})$. An $N$-derivation
$\varphi$ is called an \emph{$N$-derivation of homogeneous degree
$\alpha$} if $\varphi({\cal L}_{\beta})\subseteq {\cal
L}_{\alpha+\beta},$ for ${\bf \alpha},\beta\in G$. Set
$$\textrm{Der}^{(N)}_\alpha ({\cal L}) =\{\varphi\in \textrm{Der}^{(N)}({\cal L})|\deg \varphi={\bf \alpha}\}.
$$

\vspace{3mm} \no\textbf{Lemma 2.1}\quad {\it Let $G$ be an abelian
group. For any finitely generated $G$-graded Lie algebra
 ${\cal L}=\bigoplus\limits_{{\bf \alpha}\in G}{\cal L}_{{\bf \alpha}}$,
we have}
$$
\textrm{Der}^{(N)}({\cal L})=\bigoplus\limits_{{\bf \alpha}\in G}
\textrm{Der}^{(N)}_\alpha({\cal L}).
$$
\emph{Proof} \quad For $\alpha\in G$, let $\rho_\alpha:{\cal
L}\rightarrow{\cal L}_\alpha$ denote the canonical projection. Let
$S$ be a finite subset generating  ${\cal L}$. Set
$$Y:=S\cup\{[x_1,\cdots,x_m]|x_1,\cdots,x_m\in S,\ m=2,\cdots,N-1\},$$ $Y$ is also a finite
subset of ${\cal L}$.  Clearly, $Y$ generates  ${\cal L}$ as an
$N$-Lie algebra (i.e.,  ${\cal L}$ is the smallest subspace
containing $Y$ and stable under taking iterated brackets of the form
$[x_1,\cdots,x_N]$).
 For
$\varphi\in\textrm{Der}^{(N)}({\cal L})$, there is a finite set $
K\subseteq G$ such that
$$ Y\cup\varphi(Y)\subset \sum_{\alpha\in K}{\cal L}_\alpha.\eqno(1)$$

For $\alpha\in G$, set $\varphi_\alpha:=\sum_{\beta\in
G}\rho_{\alpha+\beta}\varphi\rho_{\beta}$. Since for
$x_{\beta_i}\in{\cal L}_{\beta_i}$ $(i=1,\cdots,N)$, we have
$$\begin{array}{ll} & \varphi_\alpha([x_{\beta_1}, \cdots,x_{\beta_N}])\\
=& \rho_{\alpha+\beta_1+\cdots+\beta_N}\varphi([x_{\beta_1}, \cdots, x_{\beta_N}])\\
 =&
 \rho_{\alpha+\beta_1+\cdots+\beta_N}(\sum_{i=1}^{N}[x_{\beta_1},\cdots,x_{\beta_{i-1}},
\varphi(x_{\beta_i}),x_{\beta_{i+1}},\cdots,
x_{\beta_N}])\\
 =& \sum_{i=1}^{N}[x_{\beta_1},\cdots,x_{\beta_{i-1}},
\rho_{\alpha+\beta_i}\varphi(x_{\beta_i}),x_{\beta_{i+1}},\cdots,
x_{\beta_N}]
\\
 =& \sum_{i=1}^{N}[x_{\beta_1},\cdots,x_{\beta_{i-1}},
\varphi_\alpha(x_{\beta_i}),x_{\beta_{i+1}},\cdots, x_{\beta_N}],
\end{array}
 $$
 which follows that $\varphi_\alpha\in \textrm{Der}^{(N)}_\alpha({\cal L})$.

Let $T:=\{\alpha-\beta|\alpha,\beta\in K\}$, then $T$ is finite. For
$y\in Y$,  we obtain
$$\begin{array}{rl}  \varphi(y) \overset{(a)}{=}& \sum_{\alpha,\beta\in
K}\rho_\alpha\varphi\rho_\beta(y)\\
=&\sum_{\alpha,\beta\in
K}\rho_{\alpha-\beta+\beta}\varphi\rho_\beta(y)
\\
=& \sum_{\beta\in K}\sum_{\gamma\in K-\beta}\rho_{\gamma+\beta}\varphi\rho_\beta(y)\\
  \overset{(b)}{=}& \sum_{\beta\in K}\sum_{\gamma\in T}\rho_{\gamma+\beta}\varphi\rho_\beta(y)\\
=& \sum_{\gamma\in T}\sum_{\beta\in K}\rho_{\gamma+\beta}\varphi\rho_\beta(y)\\
  \overset{(c)}{=}& \sum_{\gamma\in T}\sum_{\beta\in G}\rho_{\gamma+\beta}\varphi\rho_\beta(y)\\
    =& \sum_{\gamma\in T}\varphi_\gamma(y), \end{array}
 $$
where  $(a)$, $(b)$ and  $(c)$ follow from (1). This shows that the
Lie $N$-derivation $\varphi$ and $\sum_{\gamma\in T}\varphi_\gamma$
coincide on $Y$. By the construction of $Y$ and the definition of
Lie $N$-derivation, we obtain $\varphi=\sum_{\gamma\in
T}\varphi_\gamma$. \hfill $\Box$

\vspace{3mm}\no\textbf{Remark 2.2}\quad In case $N=2$,
$\textrm{Der}^{(N)}({\cal L})=\textrm{Der}({\cal L})$, the above
lemma  is contained in the proposition 1.1 of \cite{Fa}. In what
follows, we assume $N\geq 3$.

\section{Main theorem}

Let ${\cal L}$ be a finitely generated Lie algebra graded by a
nontrivial finite dimensional Cartan subalgebra $\mathcal{H}$. Let
$\mathcal{H}^*$ be the dual space of $\mathcal{H}$. For $\alpha\in
\mathcal{H}^*$,
$${\cal L}_\alpha=\{x\in{\cal L}|[h,x]=\alpha(h)x \quad \textrm{for all}
\ h\in \mathcal{H}\}$$ is the root space associated to $\alpha$. Let
$R=\{\alpha\in\mathcal{H}^*| {\cal L}_\alpha\neq\{0\}\}$, $R$ is
called the root system of $\mathcal{L}$ with respect to
$\mathcal{H}$. Thus we have
$${\cal L}=\oplus_{\alpha\in R}{\cal
L}_\alpha,\quad \textrm{and}\quad \mathcal{L}_0=\mathcal{H}.$$ Set
$R^{\times}=\{\alpha\in R| \alpha\neq 0\}$, $R_{\pm}=\{\alpha\in
R^{\times}| -\alpha\in R\}$.

\vspace{3mm}\no \textbf{Definition 3.1} \quad   \emph{We say  ${\cal
L}$ satisfying  property (P) if for $\alpha\in R_{\pm}$, every
element $x_\alpha\in\mathcal{L}_\alpha$  satisfies (P1) or (P2),
where}

\vspace{3mm}\emph{(P1)\quad
$x_{\alpha}=[y_{\alpha-\beta},y_{\beta}]$ for some
$y_{\alpha-\beta}\in\mathcal{L}_{\alpha-\beta}$ and
$y_{\beta}\in\mathcal{L}_{\beta}$, with
 $\beta\not\in\{0,\alpha\}$.}

\vspace{3mm} \emph{(P2)\quad $[x_\alpha,x_\alpha,x_{-\alpha}]\neq 0$
for some $x_{-\alpha}\in\mathcal{L}_{-\alpha}$.}

\vspace{3mm}\no \textbf{Theorem 3.2} \quad Let ${\cal L}$ be a
finitely generated Lie algebra with a nontrivial finite dimensional
Cartan subalgebra $\mathcal{H}$, and let  $N\geq3$ be a positive
integer. If $N$ is even or if ${\cal L}$  satisfying property (P),
then we have
$$\textrm{Der}^{(N)}(\mathcal{L})=\textrm{Der}(\mathcal{L}).$$
\no\emph{Proof} Assume $\mathcal{H}, R$ as above, let
$Q=\mathbb{Z}R$ be the abelian group. Clearly, $\mathcal{L}$ is a
finitely generated $Q$-graded Lie algebra. Using lemma 2.1, we have
$$
\textrm{Der}^{(N)}(\mathcal{L})=\bigoplus\limits_{{\bf \alpha}\in Q}
\textrm{Der}^{(N)}_\alpha(\mathcal{L}).
$$
Now let $\varphi\in\textrm{Der}^{(N)}_\gamma(\mathcal{L})$,
$\gamma\in Q$. In what follows we will prove $\varphi$ is a
derivation. We divide the argument into two cases.

\vspace{3mm}\textbf{Case one:}  $\gamma=0$.

For $x,y\in\mathcal{H}$,  we have $\varphi([x,y])=0=[\varphi(x),y]+
[x,\varphi(y)]$ as required. Suppose $0\neq y\in\mathcal{L}_\alpha$
for some $\alpha\in R^{\times}$, then there is
$h_\alpha\in\mathcal{H}$ such that $\alpha(h_\alpha)=1$. Taking
$\varphi$ on
$$y=[\underbrace{h_\alpha,\cdots,h_\alpha}_{N-1},y],$$ we have
$\varphi(y)=(N-1)\alpha(\varphi(h_\alpha))y+\varphi(y)$, which
implies $\alpha(\varphi(h_\alpha))=0$. Thus, for $x\in\mathcal{L}$,
taking $\varphi$ on
$$[x,y]=[x,\underbrace{h_\alpha,\cdots,h_\alpha}_{N-2},y],$$ we have
$\varphi([x,y])=[\varphi(x),y]+[x,\varphi(y)]$ as required.

 \vspace{3mm}\textbf{Case
two:}\quad  $\gamma\neq 0$.

 Since $\gamma\neq 0$, there exists $h_\gamma\in
\mathcal{H}$, such that $\gamma(h_\gamma)=1$. For $h\in
\mathcal{H}$, taking $\varphi$ on
$$0=[\underbrace{h_\gamma,\cdots,h_\gamma}_{N-1},h],$$
we have $0=[\varphi(h_\gamma),h]+\varphi(h)$, that is
$\varphi(h)=(-\textrm{ad}\varphi(h_\gamma))(h)$.
 Set $\psi=\varphi+\textrm{ad}\varphi(h_\gamma)$,  $\psi$ is an $N$-derivation and $\psi(\mathcal{H})=\{0\}$.
In what follows, we use two subcases to show
$$\psi(x_\alpha)=0,\quad \textrm{for}\  0\neq
x_\alpha\in\mathcal{L}_\alpha,\  \alpha\in R^{\times}.$$ Thus
$\varphi=-\textrm{ad}\varphi(h_\gamma)$ is a derivation.

 \vspace{3mm}\textbf{Subcase
one:}  $N$ is even or  $\gamma\neq -2\alpha$.

 Since $N$ is even or  $\gamma\neq -2\alpha$, there is $\bar{h}\in \mathcal{H}$ such that
$\alpha(\bar{h})^{N-1}\neq (\alpha+\gamma)(\bar{h})^{N-1}$.  Taking
$\psi$ on
$$\alpha(\bar{h})^{N-1}x_\alpha=[\underbrace{\bar{h},\cdots,\bar{h}}_{N-1},x_\alpha],$$
we have
$\alpha(\bar{h})^{N-1}\psi(x_\alpha)=(\alpha+\gamma)(\bar{h})^{N-1}\psi(x_\alpha)$,
which implies $\psi(x_\alpha)=0$, as required.

 \vspace{3mm}\textbf{Subcase
two:}  $N$ is odd and  $\gamma= -2\alpha$.

 Since $\alpha\neq 0$, then there
 is  $h_\alpha\in\mathcal{H}$ such that $\alpha(h_\alpha)=1$. If
 $\mathcal{L}_{-\alpha}=\{0\}$, then
we have $\psi(x_\alpha)=0$, as required. Suppose
$\mathcal{L}_{-\alpha}\neq\{0\}$, $x_\alpha\in\mathcal{L}_\alpha$.
 If $x_\alpha$ satisfies (P1), then
$x_{\alpha}=[y_{\alpha-\beta},y_{\beta}]$ for some
$y_{\alpha-\beta}\in\mathcal{L}_{\alpha-\beta}$ and
$y_{\beta}\in\mathcal{L}_{\beta}$, with
 $\beta\not\in\{0,\alpha\}$.
Using subcase one and taking $\psi$ on
$$x_{\alpha}=[\underbrace{h_\alpha,\cdots,h_\alpha}_{N-2},y_{\alpha-\beta},y_{\beta}],$$
then we have $\psi(x_{\alpha})=0$. If $x_\alpha$ satisfies (P2),
then $[x_\alpha,x_\alpha,x_{-\alpha}]\neq 0$ for some
$x_{-\alpha}\in\mathcal{L}_{-\alpha}$. Firstly, using subcase one,
we have $\psi(x_{-\alpha})=0$. Next, taking $\psi$ on
$$0=[\underbrace{h_\alpha,\cdots,h_\alpha}_{N-2},x_{-\alpha},x_{\alpha}],$$
then we have $(-2)^{N-2}[x_{-\alpha},\psi(x_{\alpha})]=0$, which
implies $[x_{-\alpha},\psi(x_{\alpha})]=0$. Finally, since $N$ is
odd, setting $h'=[x_\alpha,x_{-\alpha}]$ and taking $\psi$ on
$$\alpha(h')^{N-2}x_{\alpha}=[\underbrace{h',\cdots,h'}_{N-3},x_{\alpha},x_{-\alpha},x_{\alpha}],$$
then we have
$$\begin{array}{ll} & \alpha(h')^{N-2}\psi(x_{\alpha})\\
=&
[\underbrace{h',\cdots,h'}_{N-3},\psi(x_{\alpha}),x_{-\alpha},x_{\alpha}]\\
=& (-\alpha)(h')^{N-2}\psi(x_{\alpha})\\
=& -\alpha(h')^{N-2}\psi(x_{\alpha}).\end{array}$$ Since
$[x_\alpha,x_\alpha,x_{-\alpha}]\neq 0$, which implies
$\alpha(h')\neq 0$, we have $\psi(x_{\alpha})=0$ as required.

 Therefore, every $N$-derivations of homogeneous degree $\gamma$ are
derivations. This completes the proof. \hfill $\Box$

\vspace{3mm}\no\textbf{Remark 3.3}\quad If ${\cal L}$ is a finitely
generated Lie algebra graded by a nontrivial finite dimensional
Cartan subalgebra $\mathcal{H}$, but doesn't satisfy the property
(P), then
$\textrm{Der}^{(2N+1)}(\mathcal{L})=\textrm{Der}(\mathcal{L})$
doesn't holds in general for $N\geq 1$. See example 4.1.3 in the
following section.

\section{Applications}

\subsection{Schr\"{o}dinger-Virasoro algebra}

 The Schr\"{o}dinger-Virasoro algebra $\mathfrak{sv}$ is an infinite-dimensional Lie algebra
with $\mathbb{C}$-basis $\{L_n,M_n,Y_{n+\frac{1}{2}},C \mid n\in
\mathbb{Z}\}$ subject to the following Lie brackets:
\begin{align*}
&[L_m,L_n]=(n-m)L_{n+m}+\delta_{m+n,0}\frac{n^3-n}{12}C,\quad
[L_m,M_n]=nM_{n+m},\\
&[L_m,Y_{n+\frac{1}{2}}]=(n+\frac{1-m}{2})Y_{m+n+\frac{1}{2}},\quad [Y_{m+\frac{1}{2}},Y_{n+\frac{1}{2}}]=(n-m)M_{m+n+1},\\
&[M_m,M_n]=[M_m,Y_{n+\frac{1}{2}}]=0,\quad [\mathfrak{sv},C]=\{0\}.
\end{align*}

 It is easy to see the following facts
about $\mathfrak{sv}:$

\vskip 3mm (i) \quad $\mathfrak{sv}$ is generated by
$L_1,L_2,L_{-2},M_1$ and $Y_{-\frac{1}{2}}$.

\vskip 3mm  (ii) \quad
$\mathfrak{sv}=\oplus_{p\in\frac{1}{2}\mathbb{Z}}\mathfrak{sv}_{p}$
is a $\frac{1}{2}\mathbb{Z}$ graded Lie algebra according to the
Cartan algebra $\mathfrak{h}=\mathbb{C}L_0\oplus \mathbb{C}M_0\oplus
\mathbb{C}C$, where, $\mathfrak{sv}_{n}={\mathbb{C}}L_n\oplus
{\mathbb{C}}M_{n}$ for $n\in
\mathbb{Z}^\times(=\mathbb{Z}\backslash\{0\})$,
$\mathfrak{sv}_{n+\frac{1}{2}}=\mathbb{C} Y_{n+\frac{1}{2}}$ for
$n\in \mathbb{Z},$ and $\mathfrak{sv}_{0}=\mathfrak{h}$ .

\vskip 3mm  (iii) \quad $[L_n,L_{-n},L_n]=-2n^2L_n$, \
$[L_{2n+1},Y_{-n-\frac{1}{2}}]=(-2n-1)Y_{n+\frac{1}{2}}$, \
$[L_{-2n-1},M_{3n+1}]$ $=(3n+1)M_n$, for $n\in \mathbb{Z}^\times.$

\vskip 3mm Therefore, $\mathfrak{sv}$ is a finitely generated Lie
algebra graded by a nontrivial finite dimensional Cartan subalgebra
$\mathfrak{h}$  and satisfies property (P). Using theorem 3.2, we
have corollary 4.1.1.

\vspace{3mm}\no \textbf{Corollary 4.1.1} \quad
$\textrm{Der}^{(N)}(\mathfrak{sv})=\textrm{Der}(\mathfrak{sv}) \quad
(N\geq 3).$

\vspace{3mm}\no \textbf{Example 4.1.2} \quad
 Let $K=\textrm{span}_{\mathbb{C}}\{L_0,M_1,M_{-1}\}$ be the subalgebra of $\mathfrak{sv}$. $K$ is a finitely generated Lie
algebra with Cartan subalgebra $\mathcal{H}=\mathbb{C}L_0$, but
doesn't satisfy property (P).
 Let $\varphi:K\rightarrow K$ be a linear map such that
$$\varphi(L_0)=\varphi(M_{-1})=0, \quad  \varphi(M_{1})=M_{-1}.$$
Taking $\varphi$ on $[L_0,M_1]=M_1$, we get $\varphi$ is not a
derivation of $K$. But one can check that $\varphi$ is an
$N$-derivation of $K$ for any odd integer $N\geq3$.

\subsection{Generalied Witt algebra}

Let $\cal A$ be the Laurent polynomial ring ${\bf
C}[t_{1}^{\pm1},t_{2}^{\pm1},\cdots,t_{d}^{\pm1}]$ with commuting
variables, and let $\mathfrak{w}_d$ be the  Lie algebra of
derivations of $\cal A$. $\mathfrak{w}_d$ is called generalied Witt
algebra, and is also the Lie algebra of diffeomorphisms of torus
$T^d$ (see \cite{RSS}). When $d=1$, $\mathfrak{w}_d$ is the Witt
algebra and its universal central extension is called the Virasoro
algebra.  When $d\geq 2$, $\mathfrak{w}_d$  has no nontrivial
central extension.

 Let
$\mathbb{Z}^d=\mathbb{Z}\varepsilon_1\oplus\cdots\oplus\mathbb{Z}\varepsilon_d$.
For ${\bf n}=n_1\varepsilon_1+\cdots+n_d\varepsilon_d\in
\mathbb{Z}^d$, let $t^{\bf n}=t_1^{n_1}t_2^{n_2}\cdots t_d^{n_d}$
and let $D_j({\bf n})=t^{\bf n}t_j\frac{\partial}{\partial t_j}$.
 Then
$\mathfrak{w}_d=\textrm{span}_{\mathbb C}\{D_j({\bf n})|{\bf
n}\in\mathbb{Z}^d,\ j=1,2\cdots,d\}$ with the following Lie
structure:
$$[D_j({\bf n}),D_k({\bf
m})]=m_jD_k(\mathbf{n}+\mathbf{m})-n_kD_j(\mathbf{n}+\mathbf{m}), \
\forall \mathbf{n},\mathbf{m}\in \mathbb{Z}^d,\ j,k=1,2\cdots,d.
$$

 It is easy to see the following facts
about $\mathfrak{w}_d$:

\vskip 3mm (i) \quad $\mathfrak{w}_d$ is  generated by $\{D_i({\pm
}\varepsilon_j),D_i({\pm 2}\varepsilon_j)|i,j=1,2\cdots,d\}$.

 \vskip 3mm (ii)\quad
$\mathfrak{w}_d=\oplus_{\mathbf{n}\in
{\mathbb{Z}}^{d}}(\mathfrak{w}_d)_{\mathbf{n}}$ is a
${\mathbb{Z}}^{d}$-graded Lie algebra respects to the Cartan
subalgebra ${\cal H}=\textrm{span}_{\mathbb C}\{D_1({\bf
0}),\cdots,D_d({\bf 0})\}$, where
$(\mathfrak{w}_d)_{\mathbf{n}}=\textrm{span}_{\mathbb C}\{D_1({\bf
n}),\cdots,D_d({\bf n})\}$.

\vskip 3mm  (iii) \quad $[D_i(-(2n_i+1)\varepsilon_i),D_i({\bf
n}+(2n_i+1)\varepsilon_i)]=(5n_i+2)D_i({\bf n})$.

\vskip 3mm Therefore, $\mathfrak{w}_d$ is a finitely generated Lie
algebra graded by a nontrivial finite dimensional Cartan subalgebra
${\cal H}$  and satisfies property (P).  Using theorem 3.2, we have
corollary 4.2.

\vspace{3mm}\no \textbf{Corollary 4.2} \quad
$\textrm{Der}^{(N)}(\mathfrak{w}_d)=\textrm{Der}(\mathfrak{w}_d)$
for $N\geq 3$.

\subsection{ Kac-Moody algebra}

Recall that a matrix  $A=(a_{ij})_{i,j=1}^{n}$ is called a
generalized Cartan matrix if it satisfies the following conditions:

\vspace{2mm}\noindent (C1) \hspace{2cm} $a_{ii}=2$ for
$i=1,\cdots,n$;

\noindent (C2) \hspace{2cm} $a_{ij}$ are nonpositive integers for
$i\neq j$;

\noindent (C3) \hspace{2cm} $a_{ij}=0$ implies $a_{ji}=0$.

\vspace{2mm}Let $A=(a_{ij})_{i,j=1}^{n}$ be a generalized Cartan
matrix of rank $l$. A \emph{realization of $A$} is a triple
$(\mathfrak{h},\Pi,\Pi^{\vee})$, where $\mathfrak{h}$ is a $2n-l$
dimensional complex vector space,
$\Pi=\{\alpha_1,\cdots,\alpha_n\}\subset\mathfrak{h}^*$ and
$\Pi^{\vee}=\{\alpha_1^{\vee},\cdots,\alpha_n^{\vee}\}\subset\mathfrak{h}$
are linearly independent subsets in $\mathfrak{h}^*$ and
$\mathfrak{h}$, respectively, satisfying
$$ \alpha_j(\alpha_i^{\vee})=a_{ij} \quad \textrm{for} \quad
i,j=1,\cdots,n.
 $$

Let $\mathfrak{g}(A)$ be the Kac-Moody algebra associated to $A$.
One can see from \cite{K} that
 $\mathfrak{g}(A)$ is a complex Lie algebra
generated by $\mathfrak{h}, e_1,\cdots,e_n,f_1,\cdots,f_n$ with
following defining relations:
$$\begin{array}{ll}\  [e_i,f_j]=\delta_{ij}\alpha_i^{\vee},   &    \quad (i,j=1,\cdots,n);\\
\ [h,h']=0,   & \quad (h,h'\in \mathfrak{h});\\
\ [h,e_i]=\alpha_i(h)e_i,\ [h,f_i]=-\alpha_i(h)f_i,   &    \quad (i=1,\cdots,n; h\in\mathfrak{h}); \\
\ (\textrm{ad}e_i)^{1-a_{ij}}e_j=(\textrm{ad}f_i)^{1-a_{ij}}f_j=0,  &    \quad (i\neq j). \\
\end{array}
 $$

Let $Q=\sum_{i=1}^{n}\mathbb{Z}\alpha_i$ be the root lattice, and
$Q^+=\sum_{i=1}^{n}\mathbb{N}\alpha_i$. Then
$$\mathfrak{g}(A)=\Big(\bigoplus_{\alpha\in Q^+, \alpha\neq 0}
\mathfrak{g}_{-\alpha}\Big) \oplus  \mathfrak{h}\oplus
\Big(\bigoplus_{\alpha\in Q^+, \alpha\neq 0}
\mathfrak{g}_{\alpha}\Big)$$ is a $Q$-graded Lie algebra related to
Cartan subalgebra $\mathfrak{h}$. Here,
$\mathfrak{g}_\alpha=\{x\in\mathfrak{g}(A)|[h,x]=\alpha(h)x \
\textrm{for all} \ h\in \mathfrak{h}\}$ is a root space attached to
$\alpha$. For $\alpha>0$ (resp.  $\alpha<0$), $\mathfrak{g}_\alpha$
is the linear span of the elements of the form
$$[e_{i_1},e_{i_2},\cdots,e_{i_s}]\quad   (\textrm{resp. }[f_{i_1},
f_{i_{2}},\cdots,f_{i_s}])$$ such that
$\alpha_{i_1}+\cdots+\alpha_{i_s}=\alpha$ (resp. $=-\alpha$).
Moreover, we have
$$\begin{array}{lll}(i)&\mathfrak{g}_{\alpha_i}=\mathbb{C}e_i,\quad\mathfrak{g}_{-\alpha_i}=\mathbb{C}f_i,&
\textrm{for } i\in\{1,\cdots,n\};\\
(ii)&\mathfrak{g}_{s\alpha_i}=\{0\},& \textrm{for }i\in\{1,\cdots,n\}, |s|>1;\\
 (iii)&  [e_i,f_i,e_i]=2e_i\neq 0 ,\quad  [f_i,e_i,f_i]=2f_i\neq
0, &\textrm{for } i\in\{1,\cdots,n\};\\
(iv)&
 [e_{i_1},e_{i_2},\cdots,e_{i_s}]=[e_{i_1},[e_{i_2},\cdots,e_{i_s}]],&\\
 &[f_{i_1},f_{i_2},\cdots,f_{i_s}]=[f_{i_1},[f_{i_2},\cdots,f_{i_s}]],&
\textrm{for } i_1,\cdots,i_s\in\{1,\cdots,n\},\ s\geq
2.\end{array}$$

\vskip 3mm Let $\mathfrak{b}^{\pm}=\bigoplus_{\alpha\in Q^+}
\mathfrak{g}_{\pm\alpha}$. $\mathfrak{b}^{+}, \mathfrak{b}^{-}$ are
Borel subalgebras of $\mathfrak{g}(A)$. One can check that
$\mathfrak{g}(A)$ and $\mathfrak{b}^{\pm}$  are finitely generated
Lie algebras graded by a nontrivial finite dimensional Cartan
subalgebra $\mathfrak{h}$  and satisfy  property (P).  Using theorem
3.2, we have corollary 4.3.1.

\vspace{3mm}\no \textbf{Corollary 4.3.1} \quad Let
$A=(a_{ij})_{i,j=1}^{n}$ be a generalized Cartan matrix,
$\mathfrak{g}(A)$ and $\mathfrak{b}^{\pm}$ as above. For $N\geq 3$,
$$\textrm{Der}^{(N)}(\mathfrak{g}(A))=\textrm{Der}(\mathfrak{g}(A)),\quad
\textrm{Der}^{(N)}(\mathfrak{b}^{\pm})=\textrm{Der}(\mathfrak{b}^{\pm}).$$

\vspace{3mm}\no \textbf{Remark 4.3.2} \quad Let $\mathfrak{g}$ be a
semisimple complex Lie algebra, and $A$ be the  Cartan matrix of
$\mathfrak{g}$.  Then $\mathfrak{g}(A)=\mathfrak{g}$. Since any
derivation of $\mathfrak{g}$ is an inner derivation, by corollary
4.3.1,  any $N$-derivation of semisimple complex Lie algebra is an
inner derivation for $N\geq 3$.

\end{document}